\documentclass{elsart}
\usepackage[english]{babel}
\usepackage{amsfonts}
\usepackage{amsmath}
\usepackage{verbatim}
\usepackage{amssymb,amsbsy}
\usepackage{amscd}
\usepackage{latexsym}
\usepackage{float}
\usepackage{graphicx}
\usepackage{setspace}
\usepackage{color}
\usepackage{stmaryrd,algorithm} 
\usepackage[bordercolor=white, color=white, colorinlistoftodos]{todonotes}

\newcommand{\jump}[1]{\llbracket #1 \rrbracket }
\newcommand{\bbR}{\mathbb{R}}

\usepackage{subfigure}

\usepackage{color}
\definecolor{myorange}{rgb}{0.9568,0.4941,0.1961}
\definecolor{myred}{rgb}{0.9098,0.1294,0.2078}
\definecolor{myblue}{rgb}{0.0352,0.4981,0.6509}
\definecolor{myhyperblue}{rgb}{0.1607,0.3922,0.9}
\definecolor{mygreen}{rgb}{0.2235,0.6353,0.2588}
\definecolor{mygrey}{rgb}{0.3,0.3,0.3}

\begin{document}
% **********************************
\begin{frontmatter}

\title{Galerkin least squares finite element method for the obstacle problem}

\author[UCL]{Erik Burman,}  
\author[PH]{Peter Hansbo,}  
\author[ML]{Mats G. Larson,}  
\author[RS]{Rolf Stenberg}  
\address[UCL]{Department of Mathematics, 
University College London, London, 
UK--WC1E  6BT, UK}%; {\tt E.N.Burman@sussex.ac.uk}}
\address[PH]{Department of Mechanical Engineering, J\"onk\"oping University,
SE-55111 J\"onk\"oping, Sweden}
\address[ML]{Department of Mathematics and Mathematical Statistics, Ume{\aa} University, SE-901 87 Ume{\aa}, Sweden}
\address[RS]{Institute of Mathematics, Helsinki University of Technology, P. O. Box 1100, 02015
TKK, Finland}

\begin{abstract}
We construct a consistent multiplier free method for the finite element 
solution of the obstacle problem. The method is based on an augmented 
Lagrangian formulation in which we eliminate the multiplier by use of 
its definition in a discrete setting. We prove existence and uniqueness 
of discrete solutions and optimal order a priori error estimates for 
smooth exact solutions. Using a saturation assumption we also prove an a 
posteriori error estimate. Numerical examples show the performance of the 
method and of an adaptive algorithm for the control of the discretization 
error.
\end{abstract}


\end{frontmatter}
\maketitle

%!TEX root = /Users/mfernand/Documents/Articles/fd-fsi/paper.tex
%\HOX{Peter: Diskutera variations formulering\\
%Relation till mekaniken och penalty formulering, obstacle}
\section{Introduction}
Our aim in this paper is to design a simple consistent penalty method for 
contact problems that avoids the solution
of variational inequalities. We eliminate the need for Lagrange multipliers to
enforce the contact conditions by using its definition in a discrete setting,
following an idea of Chouly and Hild \cite{CH13b} used for elastic contact.

\subsection{The model problem}

We consider the obstacle problem of finding the displacement $u$ of a membrane
constrained to stay above an obstacle given by $\psi=\psi(x,y)$ (with $\psi \leq 0$ at $\partial\Omega$):
\begin{equation}\label{obstacle}
\begin{array}{rcl}
-\Delta u-f &\geq& 0 \mbox{ in } \Omega\subset \bbR^2 \\
u &\geq&  \psi \mbox{ in } \Omega\\
(u-\psi) (f+\Delta u) &= & 0  \mbox{ in }  \Omega\\
u &= & 0  \mbox{ on } \partial \Omega 
\end{array}
\end{equation}
where $\Omega$ is a convex polygon. 
It is well known that this problem admits a unique solutions $u \in
H^1(\Omega)$. This follows from the theory of Stampacchia applied to
the corresponding variational inequality (see for instance
\cite{HHN96}).

\subsection{The finite element method}

There exists a large body of literature treating finite
element methods for unilateral problems in general and 
obstacle problems in particular, e.g., \cite{KO88,Sch84,GleT89,Joh92,HHN96,ChNo00,Ve01,Zh01,BrCaHo07,WeWo10}. 
Discretization of \eqref{obstacle} is usually performed directly starting 
from the variational inequality or using a penalty method. The first 
approach however leads to some nontrivial
choices in the construction of the discretization spaces in order to
satisfy the nonpenetration condition and associated inf-sup conditions
and until recently it has proved difficult to
obtain optimal error estimates \cite{HR12,DH16}. The latter approach, o
n the other hand leads to the usual
consistency and conditioning problems of penalty methods. 

An alternative is to use the augmented Lagrangian method. We 
introduce 
the Lagrange multiplier $\lambda$ such that 
\begin{equation}\label{mult1}
\begin{array}{rcl}
-\Delta u +\lambda &=& f \mbox{ in } \Omega \\
u &= & 0  \mbox{ on } \partial \Omega 
\end{array}
\end{equation}
under the Kuhn-Tucker side conditions
\begin{equation}\label{mult2}
\begin{array}{rcl}
\psi -u &\leq&  0 \mbox{ in } \Omega\\
\lambda &\leq&  0 \mbox{ in } \Omega\\
(\psi-u) \,\lambda &= & 0  \mbox{ on }  \Omega .
\end{array}
\end{equation}
Using the standard trick of rewriting the Kuhn-Tucker conditions as
\begin{equation}
\lambda = -\frac{1}{\gamma}\left[ \psi - u   -\gamma\lambda\right]_+
\end{equation}
where $[x]_\pm = \pm \max(0,\pm x)$ and $\gamma\in \bbR^+$, cf., e.g., Chouly and Hild \cite{CH13b}, we can formulate the augmented Lagrangian problem of finding $(u,\lambda)$ that are stationary points to the functional
\begin{align}\nonumber
\frak{F}(u,\lambda) := {}& \frac12\int_{\Omega} \vert\nabla u\vert^2\, d\Omega + \int_{\Omega}\frac{1}{2\gamma}\left[  \psi-u-\gamma\lambda\right]_+^2\, d\Omega \\
 {}& -\frac12\int_{\Omega} \gamma \lambda^2 \, d{\Omega}
-\int_{\Omega} f u \, d{\Omega} ,
\end{align}
cf. Alart and Curnier \cite{AC91}, leading to seeking $(u,\lambda)\in H_0^1(\Omega)\times L_2(\Omega)$ such that
\begin{equation}
 \int_{\Omega} \nabla u\cdot\nabla v\, d\Omega - \int_{\Omega}\frac{1}{\gamma}\left[ \psi - u -\gamma\lambda\right]_+ v\, d\Omega =\int_{\Omega} f v \, d{\Omega}\quad\forall v\in H_0^1(\Omega)
\end{equation}
and
\begin{equation}
\int_{\Omega}\frac{1}{\gamma}\left[ \psi - u -\gamma\lambda\right]_+ \mu\, d\Omega +\int_{\Omega} \lambda \mu  \, d{\Omega}=0\quad\forall \mu\in L_2(\Omega) .
\end{equation}

%Another approach proposed by Hild and
%Renard \cite{HR10} is to use a
%stabilized Lagrange-multiplier in the spirit of Barbosa and Hughes \cite{BH92}.
%As a further development one may use the reformulation of the contact condition
%\begin{equation}\label{cond1}
%\partial_n u = - \gamma^{-1} [u - \gamma \partial_n u ]_+
%\end{equation}
% introduced by Alart and Curnier \cite{AC91} in an augmented Lagrangian
%framework.
%Using the close relationship between the Barbosa and Hughes method and
%Nitsche's method \cite{Nit70} discussed by Stenberg \cite{Sten95}, this method was
%then further developed in the elegant Nitsche-type formulation
%introduced by Chouly, Hild and Renard  \cite{CH13b,CHR15}. In these works
%optimal error estimates for the above model problem were obtained for the first time.

For our discrete method, we assume that $\{ \mathcal{T}\}_h$ is a
family of conforming shape regular meshes on $\Omega$, consisting of
triangles $\mathcal{T} = \{T\}$
and 
define $V_h$ as the space of $H^1$--conforming piecewise polynomial
functions 
on $\mathcal{T}$, satisfying the homogeneous boundary condition of
$\Gamma_D$.
\[
V_h := \{v_h \in H^1_0(\Omega): v_h\vert_T \in \mathbb{P}_k(T), \,
\forall T \in \mathcal{T} \},\quad \mbox{ for } k \ge 2.
\]
We then formally replace $\lambda$ element--wise by
$\Delta u_h + f$ to obtain a discrete minimization problem: 
seek $u_h\in V_h$ such that
\begin{equation}
u_h = \arg\min_{v\in V_h} \frak{F}_h(v)
\end{equation}
where
\begin{align}\nonumber
\frak{F}_h(v) := {}& \frac12\int_{\Omega} \vert\nabla v\vert^2\, d\Omega + \sum_{T\in\mathcal{T}}\int_{T}\frac{1}{2\gamma}\left[ \psi - v -\gamma(\Delta v + f)\right]_+^2\, d\Omega  \\
{}& -\frac12\sum_{T\in\mathcal{T}}\int_{T} \gamma (\Delta v + f)^2 \, d{\Omega}
-\int_{\Omega} f v \, d{\Omega} .\label{perturbed}
\end{align}
The
Euler--Lagrange equations corresponding to (\ref{perturbed}) take the form: Find $u_h \in V_h$ such that
\begin{equation}\label{FEM}
a(u_h,v_h) + b(u_h,\psi,f;v_h) =(f,v_h)_\Omega \quad \forall v_h
\in V_h 
\end{equation}
where $(\cdot,\cdot)_\Omega$ denotes the standard $L^2$-inner product,
$a(u_h,v_h) := (\nabla u_h,\nabla v_h)_{\Omega}$ and 
\begin{align}\nonumber
b(u_h,\psi,f;v_h):={}& \left<
 \gamma^{-1}[\psi -u_h- \gamma (\Delta u_h + f)]_+, -v_h - \gamma \Delta
 v_h\right>_{h} \\ 
 &{} -   \left<\gamma (\Delta u_h + f),  \Delta
 v_h\right>_{h} \label{stab_form0}
\end{align}
where 
\[
\left<x_h,y_h\right>_{h} := \sum_{T\in\mathcal{T}} \int_T x_h y_h~ \mbox{d} x 
\]
and, for use below,
\[
\| x_h\|_h := \left<x_h,x_h\right>_{h}^{1/2}.
\]
To simplify the notation below we introduce $P_\gamma(u_h) = 
\gamma \Delta u_h+u_h $ and
\[
b(u_h,\psi,f;v_h):=\left<
 \gamma^{-1}[\psi - \gamma f-P_\gamma(u_h)]_+, -P_\gamma(v_h) \right>_{h} -   \left<\gamma (\Delta u_h + f),  \Delta
 v_h\right>_{h}.
\]
We will also omit $\psi$ and $f$ from the
argument of $b$ below, and use the notation $\Psi:=\psi -
\gamma f$ so that
\[
b(u_h;v_h):=\left<
 \gamma^{-1}[\Psi-P_\gamma(u_h)]_+, -P_\gamma(v_h) \right>_{h} -   \left<\gamma (\Delta u_h + f),  \Delta
 v_h\right>_{h}.
\]
Note that the form $b$ can be interpreted as a nonlinear consistent
least squares penalty term for the imposition of the contact
condition. A similar method was proposed in Stenberg et al. \cite{GuStVi16} 
in the framework of variational inequalities.

We will below alternatively use the compact notation
\[
A_h(u_h,v_h):=a(u_h,v_h) + b(u_h;v_h)
\]
and the associated formulation, find $u_h\in V_h$ such that
\begin{equation}\label{CompFEM}
A_h(u_h;v_h) = (f,v_h)_\Omega,\mbox{ for all }
v_h\in V_h.
\end{equation}
%We use the original formula for the contact
%condition proposed by Alart and Curnier, $\lambda = -\gamma^{-1} [
%u-\gamma \lambda]_+$ and formally replace $\lambda$ element wise by
%$-\Delta u_h - f$ and $\mu_h$ by $-\Delta v_h$.

\subsection{Summary of main results and outline} In Section 2 we 
recall some technical 
results, in Section 3 we derive an existence result for the discrete solution 
using Brouwer's fixed point theorem and we prove uniqueness of the solution 
using monotinicity of the the nonlinearity, in Section 4 we prove an a 
priori error estimate and using a saturation assumption we also prove an 
a posteriori error estimate, finally in Section 5 we present numerical 
results confirming our theoretical 
results and illustrating the performance of an adaptive algorithm 
based on our a posteriori error estimate.

\section{Technical results}
Below we will use the notation $a \lesssim b$ for $a \leq C b$ where
$C$ is a constant independent of $h$, but not of the local mesh
geometry.

First we recall the following inverse inequality,
\begin{equation}\label{eq:inv_ineq}
\|\nabla v_h\|_{T} \leq C_i h_T^{-1} \|v_h\|_T , \qquad T \in \mathcal{T}
\end{equation}
see Thom\'ee \cite{Thomee2006}.

We will use the Scott-Zhang interpolant preserving boundary
conditions, denoted $i_h:H^1_0(\Omega) \mapsto V_h$. This operator is
$H^1$-stable, $\|i_h u\|_{H^1(\Omega)} \lesssim \|u\|_{H^1(\Omega)}$
and the following
interpolation error estimate is known to hold \cite{SZ90},
\begin{equation}\label{interpolation}
\|u-i_hu\|_{\Omega}+ h \|u-i_hu\|_{H^1(\Omega)}
+  h^2 \|
\Delta (u -i_hu)\|_{h} \lesssim  h^{k+1} |u|_{H^{k+1}(\Omega)}.
\end{equation}

The essential properties of the nonlinearity are collected in the
following lemmas.
\begin{lem}\label{lem:monotone}
Let $a,b \in \mathbb{R}$ then there holds
\[
([a]_+-[b]_+)^2 \leq ([a]_+ - [b]_+) (a - b),
\]
\[
|[a]_+-[b]_+| \leq |a-b|.
\]
\end{lem}
\begin{pf}
Developing the left hand side of the expression we have
\[
[a]_+^2 +[b]_+^2 - 2 [a]_+[b]_+ \leq  [a]_+ a + [b]_+ b -  a [b]_+
-  [a]_+b =  ([a]_+ - [b]_+) (a - b).
\]
For the proof of the second claim, this is trivially true in case both $a$
and $b$ are positive or negative. If $a$ is negative and $b$ positive then
\[
|[a]_+-[b]_+| = |b| \leq |b-a|
\]
and similarly if $b$ is negative and $a$ positive
\[
|[a]_+-[b]_+| = |a| \leq |b-a|.
\]
\end{pf}
\begin{lem}(Continuity of $b(\cdot;\cdot)$)\label{bcont}
For all $u_1,u_2,v \in V_h$, the form \eqref{stab_form0} satisfies
\begin{equation}
|b(u_1;v) - b(u_2;v)| \lesssim 
\gamma^{-1} (\| (u_1-u_2)\|_{\Omega}+ \gamma h^{-1} \|\nabla (u_1-u_2)\|_\Omega) (\|v\|_{\Omega}+ \gamma h^{-1}
\|\nabla v_h\|_\Omega).
\end{equation}
\end{lem}
\begin{pf}
\begin{align*}
b(u_1;v_h) - b(u_2;v_h) = {}& \gamma^{-1}\left<
 [ \Psi-P_\gamma(u_1)]_+- [\Psi-P_\gamma(u_2)]_+,- P_\gamma(v_h) \right>_{h}  \\
{}& -   \left<\gamma \Delta (u_1 - u_2),  \Delta
 v_h\right>_{h}.
\end{align*}
Using the second inequality of Lemma \ref{lem:monotone} we see that
the nonlinearity satisfies
\begin{multline}\label{1st_ineq}
\gamma^{-1}|\left<
 [\Psi-P_\gamma(u_1) ]_+- [\Psi-P_\gamma(u_2)]_+, -P_\gamma(v_h)\right>_{h}| \\
\leq \gamma^{-1} \| P_\gamma(u_1) + \Psi- P_\gamma(u_2) -\Psi\|_{C,f} \|P_\gamma(v_h)\|_{h}.
\end{multline}
By the inverse inequality \eqref{eq:inv_ineq} we have 
\begin{equation}\label{2nd_ineq}
\left<\gamma \Delta (u_1 - u_2),  \Delta
 v_h\right>_{h} \leq C_i^2 \gamma h^{-2} \|\nabla (u_1 - u_2)\|_\Omega
\|\nabla v_h\|_\Omega
\end{equation}
and
\begin{equation}\label{3rd_ineq}
\|P_\gamma(v_h)\|_{h} \leq \|v_h\|_\Omega + C_i \gamma h^{-1} \|\nabla v_h\|_\Omega.
\end{equation}
Collecting \eqref{1st_ineq},\eqref{2nd_ineq} and \eqref{3rd_ineq} we have
\begin{align*}
|b(u_1;v_h) - b(u_2;v_h) | \leq {}& (\|u_1 - u_2\|_\Omega + C_i \gamma
h^{-1} \|\nabla (u_1 - u_2)\|_\Omega) (\|v_h\|_\Omega + C_i \gamma
h^{-1} \|\nabla v_h\|_\Omega)\\
{}& + C_i^2 \gamma h^{-2} \|\nabla (u_1 - u_2)\|_\Omega
\|\nabla v_h\|_\Omega
\end{align*}
and the claim follows.
\end{pf}
\section{Existence of unique discrete solution}
In the previous works on Nitsche's method existence and uniqueness has
been proven by using the monotonicity and hemi-continuity of the
operator. Here we propose a different approach where we use the
Brouwer's fixed point theorem to establish existence and the 
monotonicity of the nonlinearity for uniqueness. We start by showing
some positivity results and a priori bounds. Since we are interested
in existence and uniqueness for a fixed mesh parameter $h$, we do not
require that the bounds in this section are uniform in $h$.
\begin{lem}\label{lem:monotonicity}
Let $u_1,u_2 \in V_h$ and assume that
\begin{equation}\label{pos_cond}
\gamma < C_i^{-2} h^2/2
\end{equation}
 then there holds
\begin{multline*}
\frac{\alpha}{2} \|u_1-u_2\|_{H^1(\Omega)}^2 +\gamma^{-1} \|[\Psi-P_\gamma(u_1) ]_+- [\Psi-P_\gamma(u_2)]_+\|^2_{h}\\
\leq A_h(u_1;u_1-u_2)-A_h(u_2;u_1-u_2)
\end{multline*}
and
\[
\frac{\alpha}{4} \|u_1\|_{H^1(\Omega)}^2 
\leq A_h(u_1;u_1) + C \alpha^{-1} \gamma^{-2} \|[\Psi]_+\|^2_\Omega.
\]
\end{lem}
\begin{pf}
First observe consider the form $b(\cdot;\cdot)$,
\begin{multline*}
b(u_1;v_h) - b(u_2;v_h) \\=\gamma^{-1}\left<
 [ \Psi-P_\gamma(u_1)]_+- [\Psi-P_\gamma(u_2)]_+, -v_h-\gamma \Delta
 v_h + \Psi - \Psi\right>_{h}  \\
-   \left<\gamma \Delta (u_1 - u_2),  \Delta
 v_h\right>_{h}.
\end{multline*}
Using the monotonicity of Lemma \ref{lem:monotone}  we may write
\begin{align*}
b(u_1;u_1-u_2) - b(u_2;u_1-u_2) 
\ge {}& \gamma^{-1}\|[ \Psi-P_\gamma(u_1) ]_+- [\Psi - P_\gamma(u_2) ]_+\|^2_{h}\\ &{} -   \gamma \|\Delta (u_1 - u_2)\|^2_{h}.
\end{align*}
Observe that using an inverse inequality \eqref{eq:inv_ineq} we have $$\gamma \|\Delta (u_1 - u_2)\|^2_{h} \leq C_i \gamma
h^{-2} \|\nabla(u_1 - u_2) \|^2_\Omega$$. 
We may then write
\begin{multline*}
(1 -  C^2_i 
h^{-2}\gamma) \|\nabla (u_1-u_2)\|_{\Omega}^2 + \gamma^{-1}\|[\Psi-P_\gamma(u_1)
]_+- [\Psi-P_\gamma(u_2) ]_+\|^2_{h}\\
\leq A_h(u_1;u_1-u_2)-A_h(u_2;u_1-u_2)
\end{multline*}
It follows that choosing
$\gamma < C_i^{-2} h^2 2$ and applying the Poincar\'e inequality
\begin{equation}\label{eq:poincare}
\alpha^{\frac12} \|u\|_{H^1(\Omega)} \leq \|\nabla u\|_\Omega, \quad
\forall u \in H^1_0(\Omega)
\end{equation}
 there holds
\begin{multline*}
\frac{\alpha}{2} \|u_1-u_2\|_{H^1(\Omega)}^2 + \gamma^{-1}\|[\Psi-P_\gamma(u_1) ]_+- [\Psi-P_\gamma(u_2) ]_+\|^2_{h}\\
\leq A_h(u_1;u_1-u_2)- A_h(u_2;u_1-u_2).
\end{multline*}
The second inequality follows by taking $u_2 = 0$ above and noting that
then
\begin{align*}
&\frac{\alpha}{2} \|u_1\|_{H^1(\Omega)}^2 +\gamma^{-1} \|[\Psi-P_\gamma(u_1)]_+- [\Psi]_+\|^2_{h}
\\
&\qquad \leq A_h(u_1;u_1)- \gamma^{-1} \left< [\Psi]_+,- P_\gamma(
 u_1) \right>_{h} \\
&\qquad \leq A_h(u_1;u_1) + (1+C_i \gamma h^{-1})\alpha^{-\frac12} \gamma^{-1}
\|[\Psi]_+\|_\Omega  \alpha^{\frac12} \|u_1\|_{H^1(\Omega)}
\end{align*}
where we used \eqref{3rd_ineq} in the last step. Considering the
condition on $\gamma$ and using an arithmetic-geometric inequality we
may conclude.
\end{pf}
\begin{prop}
Assume that  $\gamma$ saisfy \eqref{pos_cond}.
Then formulation \eqref{CompFEM} using the contact operator 
\eqref{stab_form0}, admits a unique solution.
\end{prop}
\begin{pf}
The uniqueness is an immediate consequence of Lemma \ref{lem:monotonicity}. If $u_1$ and $u_2$ both are solution to \eqref{CompFEM},
then 
\[
A_h(u_1;u_1-u_2)-A_h(u_2;u_1-u_2) = (f,u_1-u_2)_\Omega
-(f,u_1-u_2)_\Omega = 0
\]
and we conclude that $\|u_1-u_2\|_{H^1(\Omega)}=0$ and hence
$u_1\equiv u_2$.
Let $N_V$ denote the number of
degrees of freedom of $V_h$.

Consider the mapping $G:\mathbb{R}^{N_V} \mapsto \mathbb{R}^{N_V}$ defined by
\[
(G(U),V)_{\mathbb{R}^{N_V}}  := A_h(u_h,v_h)-
(f,v_h)_\Omega,
\]
where $U = \{u_i\}_{i=1}^{N_V}$, $V=\{v_i\}_{i=1}^{N_V}$,
where $\{u_i\}$ and $\{v_i\}$ denotes the vectors of unknown
associated to the basis functions of $V_h$.

By the second claim of Lemma \ref{lem:monotonicity},
there holds
\[
\frac{\alpha}{4} \|u_h\|^2_{H^1(\Omega)} - C\alpha^{-1} \gamma^{-2}  \|[\Psi]_+\|_\Omega^2 -  (f,u_h)_\Omega \leq
A_h(u_h,u_h) - (f,u_h)_\Omega = (G(U),U)_{\mathbb{R}^{N_V}}. 
\]
Since 
\[
\frac{\alpha}{4} \|u_h\|^2_{H^1(\Omega)}  - (f,u_h)_\Omega \ge
\frac{\alpha}{8} \| u_h\|^2_{H^1(\Omega)} -  C\frac{1}{\alpha} \|f\|_\Omega^2
\]
we have that for any fixed $h$ the following positivity holds for $U$ sufficiently large 
\[
0 < \frac{\alpha}{8} \|u_h\|^2_{H^1(\Omega)}  - \frac{C}{\alpha}
(\gamma^{-2}\|[\Psi]_+\|_\Omega^2+\|f\|_\Omega^2) \leq (G(U),U)_{\mathbb{R}^{N_V}}.
\]
Assume that this positivity holds whenever $|U|\ge q \in
\mathbb{R}_+$. Denote by $B_q$ the (closed) ball in
$\mathbb{R}^{N_V}$ with radius $q$
and assume that there is no $U \in B_q$ such that $G(U)=0$. Define the
function
\[
\phi(U) = - q G(U)/|G(U)|.
\]
Then $\phi:B_{q} \mapsto B_{q}$, $\phi$ is
continuous by Lemma \ref{bcont} and the assumption that $|G(U)|>0$ for
all $U \in B_q$. Hence there
exists a fixed point $X \in B_q$ such that 
\[
X = \phi(X).
\]
It follows that 
\[
|X|^2 = -q (G(X),X)/|G(X)|, 
\]
but since $|X|=q$, by assumption $(G(X), X)>0$, which leads to a
contradiction, since $|X|>0$.
\end{pf}
\section{Error estimates}

\begin{thm}(A priori error estimate)\label{thm:apriori}
Assume that $u \in H^1_0(\Omega)$ with $\Delta u \in L^2(\Omega)$ is
the solution of \eqref{obstacle} and that $u_h$ is the solution to \eqref{FEM}
with \eqref{stab_form0} and $0<\gamma =\gamma_0 h^{2}$,
where $\gamma_0 \in \mathbb{R}$, $\gamma_0 < C_i^{-2}/2$.
then there holds for all $v_h \in
V_h $
\begin{multline}\label{best_approximation}
\alpha\|u-u_h\|_{H^1(\Omega)}^2 +\gamma^{-1} \|  [\Psi-P_{\gamma}(u_h)]_+- [\Psi-P_{\gamma}(u)]_+\|^2_{h}  \\
\lesssim \frac{1}{ \alpha} \|u-v_h\|_{H^1(\Omega)}^2+  
\|\gamma^{-\frac12}(u-v_h)\|_{\Omega}
+  \|\gamma^{\frac12}
\Delta (u -  v_h)\|_{h}^2.
\end{multline}
If in addition $u \in H^{k+1}(\Omega)$ then there holds
\begin{equation}\label{optimal_convergence}
\alpha\|u-u_h\|_{H^1(\Omega)} +\gamma^{-1/2} \|
[\Psi-P_{\gamma}(u_h)]_+- [\Psi-P_{\gamma}(u)]_+\|_{h} \lesssim h^k |u|_{H^{k+1}(\Omega)}.
\end{equation}
\end{thm}
\begin{pf}
Using the definition of $a(\cdot,\cdot)$ we may write
\begin{align*}
\|\nabla (u-u_h)\|_{\Omega}^2 
&\leq a(u-u_h,u-u_h) 
\\
&= a(u-u_h,u - v_h) +a(u-u_h,v_h-u_h) 
\\
&\leq \frac{\alpha}{4} \|u-u_h\|_{H^1(\Omega)}^2
+ \frac{1}{\alpha} \|u-v_h\|_{H^1(\Omega)}^2 +a(u-u_h,v_h-u_h).
\end{align*}
Observe that
\begin{align}\nonumber
a(u,v_h-u_h) = {}& \left<-\Delta u - f + f, v_h-u_h\right>_{V',V} 
\\
=  {}& \left<\gamma^{-1} [\Psi-P_{\gamma}(u)]_+,
 (v_h-u_h) \right>_{V',V} + (f, v_h-u_h)_{\Omega}.
\end{align}
If $ [\Psi-P_{\gamma}(u)]_+ \in L^2(\Omega)$ we may also write
\[
\left<\Delta u + f, \gamma \Delta(v_h-u_h)\right>_h+\left<\gamma^{-1} [\Psi-P_{\gamma}(u)]_+,
 \gamma \Delta(v_h-u_h)\right>_h = 0.
\]
It follows that
\begin{align}\nonumber
a(u,v_h-u_h) 
= {}& (f, v_h-u_h)_{\Omega} - \left<\gamma^{-1} [\Psi-P_{\gamma}(u)]_+,- P_{\gamma}(v_h-u_h)\right>_h \\ \nonumber &{} + \left<\Delta u + f, \gamma \Delta (v_h-u_h)\right>_h\\
= {}& (f,v_h-u_h)_{\Omega} - b(u;v_h-u_h).
\end{align}
As a consequence we have the following property reminiscent of Galerkin
orthogonality,
\begin{align}\nonumber
&a(u-u_h,v_h-u_h) 
\\ \nonumber
&\qquad =  b(u_h;v_h-u_h)- b(u;v_h-u_h)
\\ \label{eq:GO_pert}
&\qquad =\left<\gamma^{-1}
  [\Psi-P_{\gamma}(u_h)]_+-\gamma^{-1} [\Psi-P_{\gamma}(u)]_+,
 -P_{\gamma}(v_h-u_h) \right>_{h} 
 \\ \nonumber
&\qquad \qquad - \gamma \left<\Delta (u_h - u), \Delta (v_h - u_h) \right>_{h}
\end{align}
First observe that
\begin{align*}
&\gamma\left<\Delta (u_h - u), \Delta (v_h - u_h) \right>_{h}
\\
&\qquad 
\leq 
\|\gamma^{\frac12}(\Delta u_h - \Delta v_h )\|^2_{h} + \|\gamma^{\frac12}
(\Delta v_h - \Delta u )\|_{h}\|\gamma^{\frac12} (\Delta v_h - \Delta u_h )\|_{h} 
\\
&\qquad 
\leq 
\frac12 C^2_i h^{-2} \gamma \|\nabla (u_h - v_h)\|^2_\Omega +  \|\gamma^{\frac12}
(\Delta v_h - \Delta u )\|_{h}^2
\\
&\qquad 
\leq 
C^2_i h^{-2} \gamma \|\nabla (u_h - u)\|^2_\Omega + C^2_i h^{-2}
\gamma \|\nabla (v_h - u)\|^2_\Omega + \|\gamma^{\frac12}
(\Delta v_h - \Delta u )\|_{h}^2.
\end{align*}
Considering the first term in the right hand side of equation
\eqref{eq:GO_pert} we may write
\begin{align*}
&\left<\gamma^{-1}
  [\Psi-P_{\gamma}(u_h)]_+-\gamma^{-1} [\Psi-P_{\gamma}(u)]_+,
  -P_{\gamma}(v_h-u_h) \right>_{h} \\
\\
&\qquad = \underbrace{\left<\gamma^{-1}
  [\Psi-P_{\gamma}(u_h)]_+-\gamma^{-1}
  [\Psi-P_{\gamma}(u)]_+,
  -P_{\gamma}(v_h-u) \right>_{h} }_{I}
\\
&\qquad \qquad + \underbrace{\left<\gamma^{-1}
  [\Psi-P_{\gamma}(u_h)]_+-\gamma^{-1}
  [\Psi-P_{\gamma}(u)]_+,
  -P_{\gamma}(u-u_h)  \right>_{h} }_{II}
\\
&\qquad = I+II.
\end{align*}
The term $I$ may be bounded using the Cauchy-Schwarz inequality
followed by the arithmetic geometric inequality
\[
I \leq \epsilon \gamma^{-1}\| [\Psi-P_{\gamma}(u_h)]_+- [\Psi-P_{\gamma}(u)]_+\|_{h}^2 + \frac{1}{4 \epsilon}
\|\gamma^{-\frac12} P_{\gamma}(v_h-u)\|_{h}^2.
\]
For the term $II$ we use the monotonicity property $([a]_+-[b]_+)(b-a)
\leq -([a]_+-[b]_+)^2$, with $a=\Psi-P_{\gamma}(u_h)$ and
$b=\Psi-P_{\gamma}(u)$ so that
\[
([a]_+-[b]_+)(b-a) = ([\Psi-P_{\gamma}(u_h)]_+-\gamma^{-1}
  [\Psi-P_{\gamma}(u)]_+)(
  \Psi -P_{\gamma}(u) - \Psi + P_\gamma(u_h))
\] to deduce that
\[
II \leq -\gamma^{-1} \|  [\Psi-P_{\gamma}(u_h)]_+- [\Psi-P_{\gamma}(u)]_+\|^2_{h}.
\]
Collecting the above bounds and using the Poincar\'e inequality
\eqref{eq:poincare} we find,
\begin{multline}
\alpha \left(\frac34 - C^2_i h^{-2} \gamma \right) \|u-u_h\|_{H^1(\Omega)}^2 \\+ \left(1-\epsilon \right) \gamma^{-1} \|  [\Psi-P_{\gamma}(u_h)]_+- [\Psi-P_{\gamma}(u)]_+\|^2_{h}
+ \left(1-\epsilon \right) \|\gamma^{\frac12}\Delta (u - u_h) \|^2_{h}  \\
\leq \frac{1}{\alpha} \|u-v_h\|_{H^1(\Omega)}^2+ \frac{1}{4 \epsilon}
\|\gamma^{-\frac12} P_{\gamma}(u-v_h)\|_{h}^2 
+ \|\gamma^{\frac12} \Delta (u-v_h) \|^2_{h}  
\end{multline}
Fixing $\epsilon=\frac12$, and fixing $\gamma$ sufficiently small so that $C_i^2h^{-2}
\gamma \leq \alpha/4$ then there holds
\begin{multline}
\alpha\|u-u_h\|_{H^1(\Omega)}^2 \gamma^{-1} \|  [\Psi-P_{\gamma}(u_h)]_+- [\Psi-P_{\gamma}(u)]_+\|^2_{h}  \\
\lesssim \frac{1}{ \alpha} \|u-v_h\|_{H^1(\Omega)}^2+  
\|\gamma^{-\frac12}(u-v_h)\|_{h}^2
+  \|\gamma^{\frac12}
\Delta (u -v_h)\|_{h}^2.
\end{multline}
This concludes the proof of \eqref{best_approximation}. The error
estimate \eqref{optimal_convergence} then follows by choosing $v_h$ to
be the interpolant, $i_h u$, applying the approximation error
estimate \eqref{interpolation} on the form
\begin{align*}
&\|u-i_hu\|_{H^1(\Omega)}+  
\|\gamma^{-\frac12}(u-i_hu)\|_{h}
+  \|\gamma^{\frac12}
\Delta (u -i_hu)\|_{h} 
\\
&\qquad \lesssim ( h^{k}+ \gamma^{-1/2} h^{k+1}+
\gamma^{1/2} h^{k-1}) |u|_{H^{k+1}(\Omega)}
\end{align*}
and using the bound on $\gamma$.
\end{pf}
{\bf Assumption:} \emph{(Saturation) We assume that there exists a constant 
$C_s$ such that}
\begin{equation}\label{saturation}
\|\Delta (u - u_h)\|_{h} \leq C_s h^{-1} \|\nabla (u - u_h)\|_\Omega.
\end{equation}
\begin{thm}(A posteriori error estimate)
Assume that $u \in H^1_0(\Omega)$ with $\Delta u \in L^2(\Omega)$ is
the solution of \eqref{obstacle} and $u_h$ the solution of \eqref{FEM}
satisfying \eqref{saturation} and with the parameter $\gamma$ satisfying $\gamma \leq \tfrac12 C_s^{-2} h^2$, then 
\begin{equation}
\alpha\|u-u_h\|_{H^1(\Omega)}+\gamma^{-1} \|  [\Psi-P_{\gamma}(u_h)]_+- [\Psi-P_{\gamma}(u)]_+\|_{h} \lesssim E(h,\gamma,u_h,f),
\end{equation}
where
\[
E(h,\gamma,u_h,f):= h \|f+\Delta u_h +\gamma^{-1}
[\Psi-P_{\gamma}(u_h)]_+\|_{h} +  \|h^{\frac12} \jump{\partial_n u_h}\|_{\mathcal{F}}.
\]
\end{thm}
\begin{pf}
Let $e=u-u_h$ then, under the assumption \eqref{saturation} and using
\eqref{lem:monotonicity} we may write
\begin{align}\nonumber
&(1 - C_s^2 h^{-2}\gamma )\|\nabla e\|_\Omega^2 + \gamma^{-1} \|  [\Psi-P_{\gamma}(u)]_+- [\Psi-P_{\gamma}(u_h)]_+\|^2_{h}
\\ \label{eq:sat_bound1}
&\qquad \leq (\nabla (u-u_h),\nabla e)_\Omega 
\\ \nonumber
&\qquad \qquad + \gamma^{-1} \left< [\Psi-P_{\gamma}(u)]_+- [\Psi-P_{\gamma}(u_h)]_+,- P_{\gamma}(e) \right>_{h}
\\ \nonumber
&\qquad \qquad -\gamma \left< \Delta (u - u_h), \Delta e\right>_{h}
\end{align}
Now, using similar arguments as in Theorem \ref{thm:apriori} we deduce
\[
 (\nabla u,\nabla e)_\Omega - \gamma^{-1}  \left<
   [\Psi-P_{\gamma}(u)]_+, P_{\gamma}(e) \right>_{h}-\gamma \left<
\Delta u, \Delta e \right>_{h} = \left<f,P_{\gamma}(e) \right>_{h}.
\]
Therefore the bound \eqref{eq:sat_bound1} may be written, under the assumption
$C_s^2 h^{-2} \gamma \leq 1/2$,
\begin{align}\nonumber
&\frac12 \|\nabla e\|_\Omega^2 + \gamma^{-1} \|  [\Psi-P_{\gamma}(u)]_+- [\Psi-P_{\gamma}(u_h)]_+\|^2_{h}
\\ \label{eq:sat_bound2}
&\qquad \leq 
\left<f,P_{\gamma}(e) \right>_{h} - (\nabla u_h,\nabla e)_\Omega 
\\ \nonumber
&\qquad\qquad  + \gamma^{-1}
\left< [\Psi-P_{\gamma}(u_h)]_+, P_{\gamma}(e) \right>_{h}+\gamma\left<
    \Delta u_h , \Delta e \right>_{h} 
\\ \nonumber
&\qquad = \left<f+\Delta u_h +\gamma^{-1} [\Psi-P_{\gamma}(u_h)]_+,(I +\gamma \Delta) (e - i_h e) \right>_h 
\\ \nonumber
&\qquad \qquad + \left<\jump{\partial_n u_h},e - i_h e
\right>_{\mathcal{F}}
\end{align}
Using that $\gamma \|\Delta i_h e\|_h \lesssim \gamma h^{-1} \|\nabla
i_h e\|_\Omega$, the $H^1$ stability of $i_h$ and the definition
of $\gamma$ we obtain
\begin{align*}
&\frac12 \|\nabla e\|_\Omega^2 + \gamma^{-1} \| [\Psi-P_{\gamma}(u)]_+- [\Psi-P_{\gamma}(u_h)]_+\|^2_{h}
\\
&\qquad 
\leq 
C (h +\gamma h^{-1}+\gamma^{\frac12})( \|f+\Delta u_h -\gamma^{-1} [\Psi-P_{\gamma}(u_h)]_+\|_{h}
\\
&\qquad \qquad \qquad \times (\|\nabla e\|_\Omega + \gamma^{\frac12} \|\Delta e\|_h)
\\
&\qquad \qquad+ C \|h^{\frac12} \jump{\partial_n u_h}\|_{\mathcal{F}} \|\nabla
e\|_\Omega  
\\
&\qquad 
\leq 
C (h +\gamma h^{-1} +\gamma^{\frac12})( \|f+\Delta u_h
-\gamma^{-1} [\Psi-P_{\gamma}(u_h)\|_{h} + \|h^{\frac12}
\jump{\partial_n u_h}\|_{\mathcal{F}})
\\
&\qquad \qquad \qquad \times(1+
C_s \gamma^{\frac12} h^{-1}) \|\nabla
e\|_\Omega
\end{align*}

Once again applying the assumption on $\gamma$ we obtain the bound
\begin{multline*}
\|\nabla e\|_\Omega + \gamma^{-1} \| [\Psi-P_{\gamma}(u)]_+- [\Psi-P_{\gamma}(u_h)]_+\|_{h} \\
\leq  C h \|f+\Delta u_h -\gamma^{-1}
[\Psi-P_{\gamma}(u_h)]_+\|_{h} +  C\|h^{\frac12} \jump{\partial_n u_h}\|_{\mathcal{F}}.
\end{multline*}
We conclude the proof using the Poincar\'e inequality \eqref{eq:poincare}.
\end{pf}
\begin{rem}
This a posteriori error estimate has the disadvantage of the
saturation assumption and also that the parameter $\gamma$ depends on
the constant in the saturation assumption. However as we shall see
below it appears to give a very good representation of the $H^1$-error
and can be used to drive adaptive refinement.
\end{rem}

\section{Numerical examples}

\subsection{A smooth rotational symmetric exact solution}
This example, from \cite{NoSiVe03}, is posed on the square $\Omega=(-1,1)\times(-1,1)$ with $\psi=0$ and
\[
f=\left\{\begin{array}{c}
    -8 r_0^2 (1-(r^2-r_0^2)) \quad \text{if $r\leq r_0$},\\
  -8 (r^2+ (r^2-r_0^2)) \quad \text{if $r> r_0$},
\end{array}\right.
\]
where $r=\sqrt{x^2+y^2}$ and $r_0=1/4$, and
with Dirichlet boundary conditions taken from the corresponding exact solution 
\[
u = [r^2-r_0^2]_+^2 .
\]

\begin{figure}[h]
\begin{center}
\includegraphics[height=7cm]{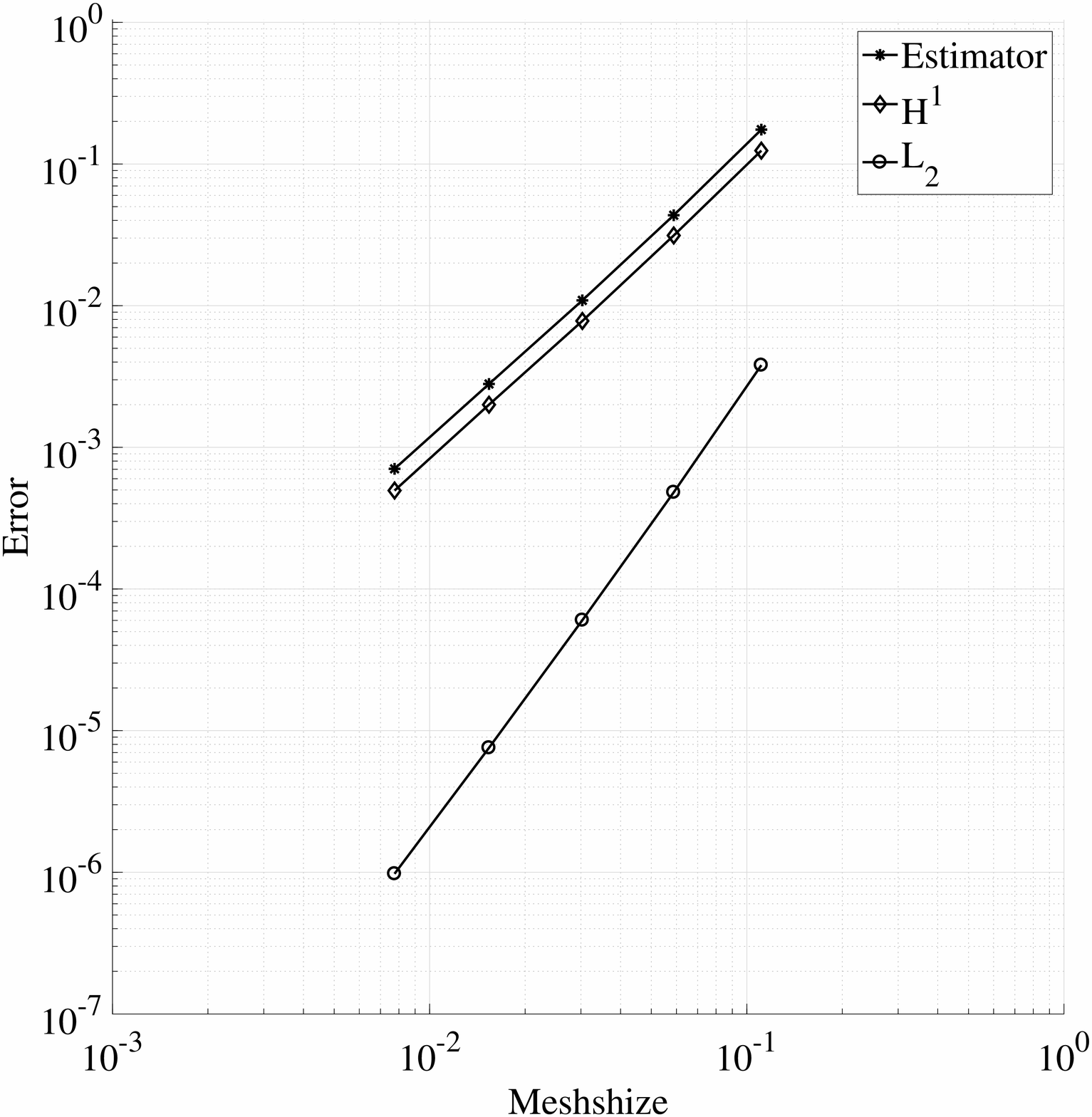}
\end{center}
\caption{Convergence for the smooth case.\label{fig:smooth}}
\end{figure}

\begin{figure}[h]
\begin{center}
\includegraphics[height=8cm]{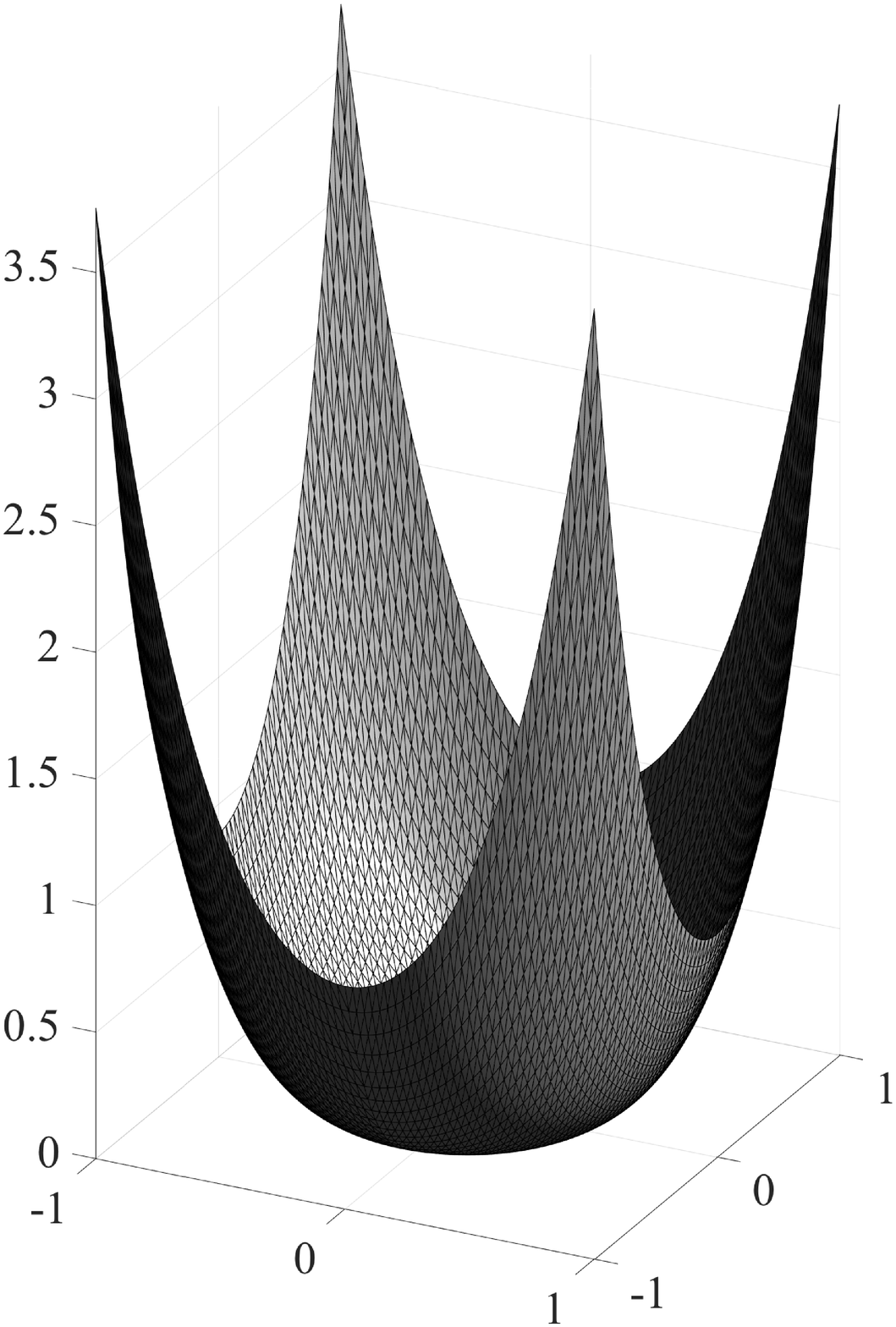}
\end{center}
\caption{Elevation of the discrete solution, smooth case.\label{fig:elevetionsmooth}}
\end{figure}

We choose $\gamma=\gamma_0 h^2$ with $\gamma_0 = 1/100$ and show the convergence in the $L_2$-- and $H^1$--norms in Figure \ref{fig:smooth} together with the 
error indicator. We remark that the smoothness of the solution precludes mesh zoning for this example, but that the indicator has the same asymptotic behaviour as the $H^1$ error. An elevation of the computed solution on one of the meshes in a sequence is given in Fig. \ref{fig:elevetionsmooth}. We note the optimal convergence of $O(h^3)$ in $L_2$ and $O(h^2)$ in $H^1$. Here we use $h=1/\sqrt{\text{NNO}}$, wehre NNO denotes the number of nodes in a uniformly refined mesh.

\subsection{A non-smooth exact solution}

This example, from \cite{BrCaHo07}, is posed on the L-shaped domain $\Omega= (-2, 2)\times (-2, 2) \setminus [0, 2)\times (-2, 0]$ with $\psi=0$ and
\[
f(r,\varphi) = -r^{2/3}\sin{(2\varphi/3)}(\gamma'(r)/r+\gamma''(r))-\frac{4}{3}r^{-1/3}\gamma'(r)\sin(2\varphi/3)-\gamma_2(r)
\]
(note the sign error in \cite{BrCaHo07}), where, with $\hat{r}=2(r-1/4)$,
\[
\gamma_1(r)=\left\{\begin{array}{ll}1,& \hat{r} < 0\\
-6\hat{r}^5+15\hat{r}^4-10\hat{r}^3+1, & 0\leq\hat{r}<1\\
0, & \hat{r}\geq 1,\end{array}\right.
\]
\[
\gamma_2(r)=\left\{\begin{array}{ll}0,& r\leq 5/4 ,\\
1 & \text{elsewhere.}\end{array}\right.
\]
with Dirichlet boundary conditions taken from the corresponding exact solution 
\[
u(r,\varphi)=r^{2/3}\gamma_1(r)\sin(2\varphi/3)
\]
which belongs to $H^{5/3-\varepsilon}(\Omega)$ for arbitrary $\varepsilon > 0$.

\begin{figure}[h]
\begin{center}
\includegraphics[height=7cm]{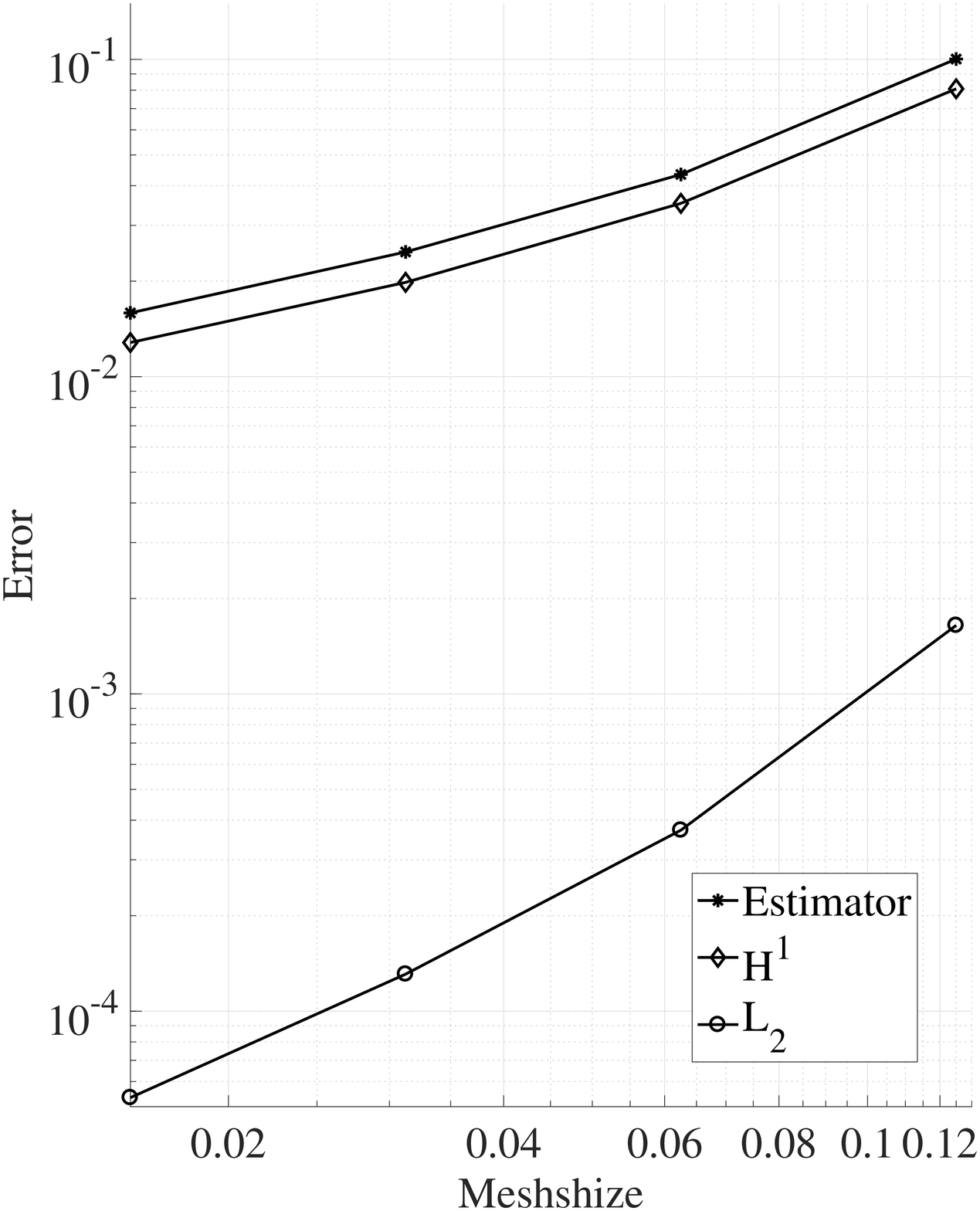}
\end{center}
\caption{Convergence for the nonsmooth case.\label{fig:nonsmootherror}}
\end{figure}

\begin{figure}[h]
\begin{center}
\includegraphics[height=5cm]{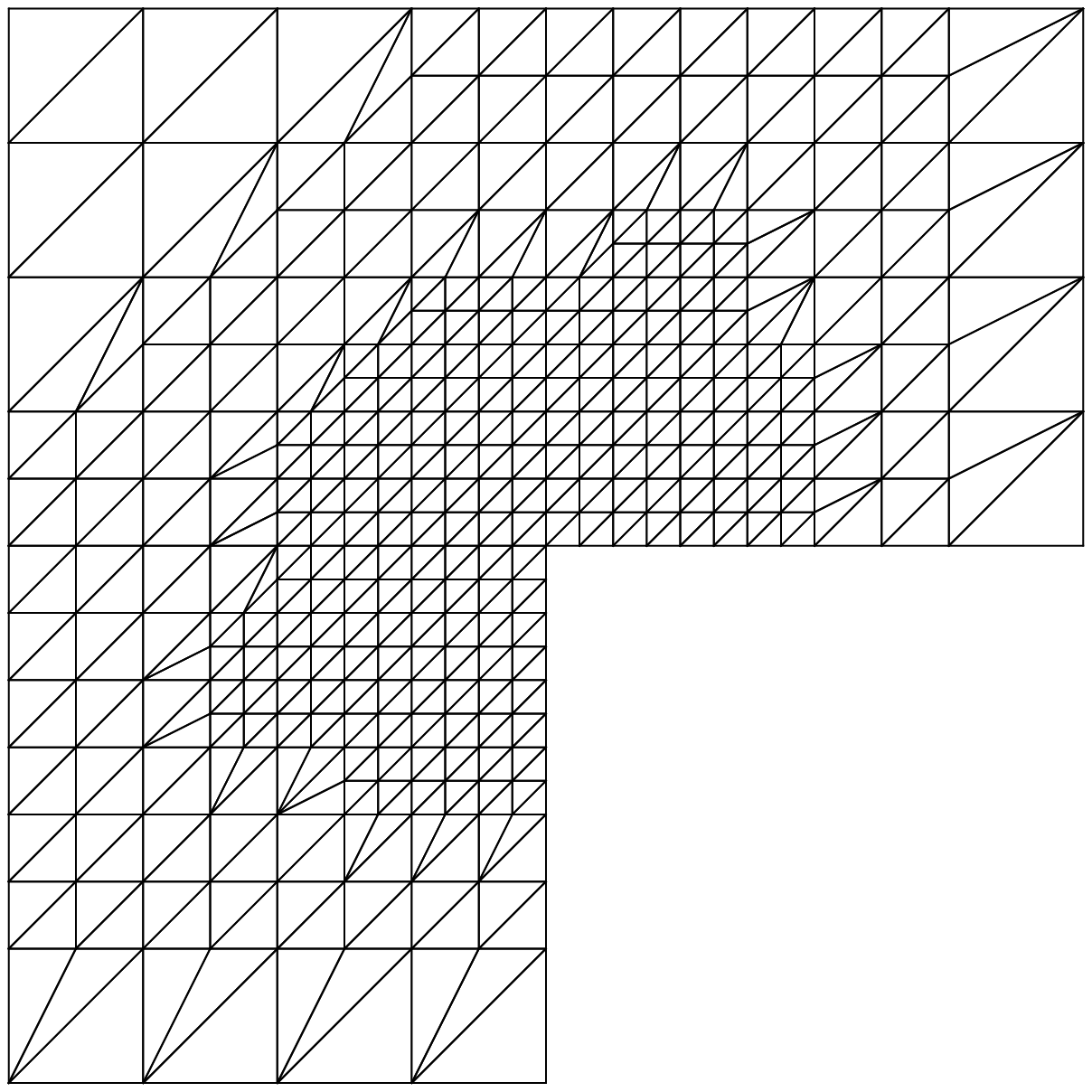}\includegraphics[height=5cm]{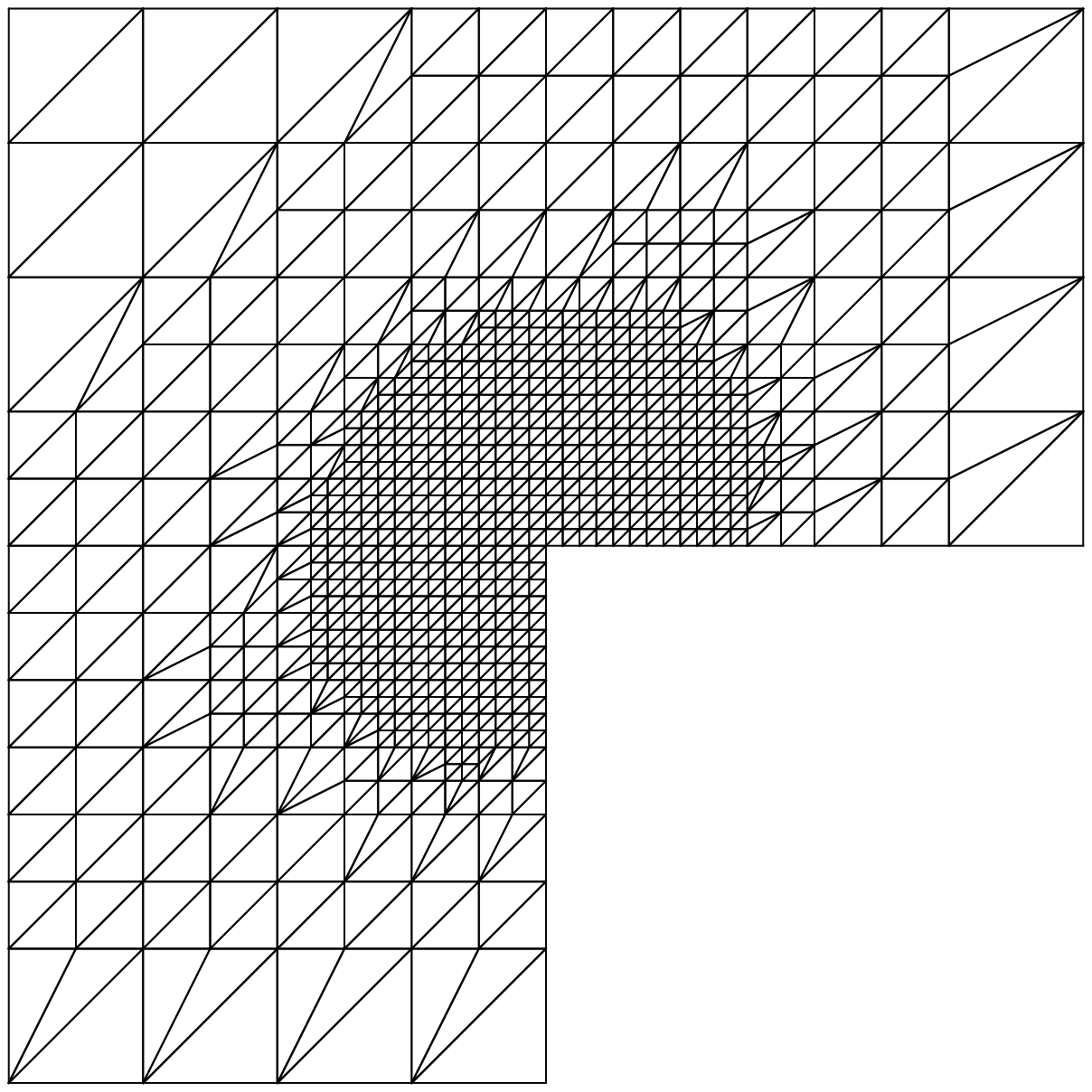}\newline\includegraphics[height=5cm]{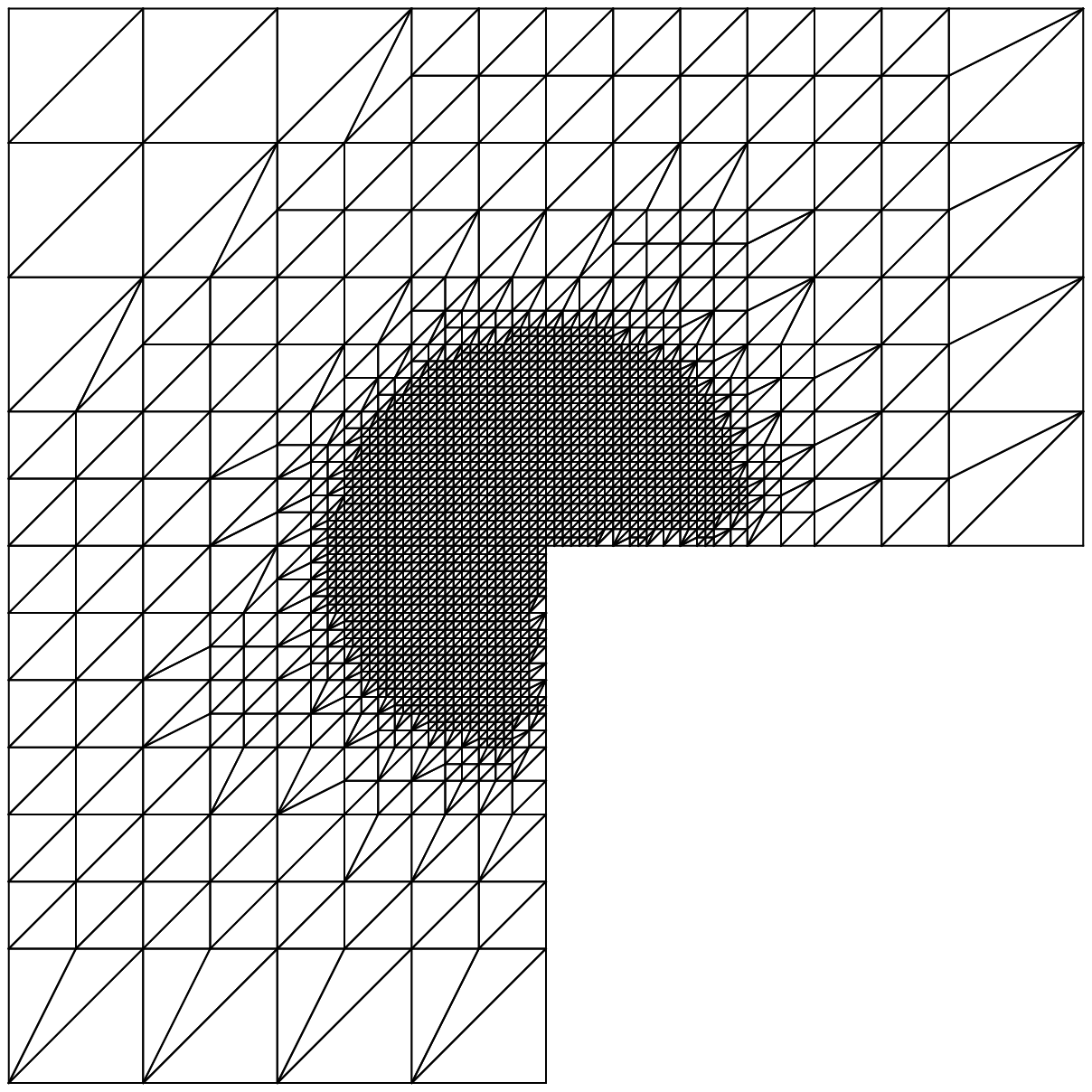}\includegraphics[height=5cm]{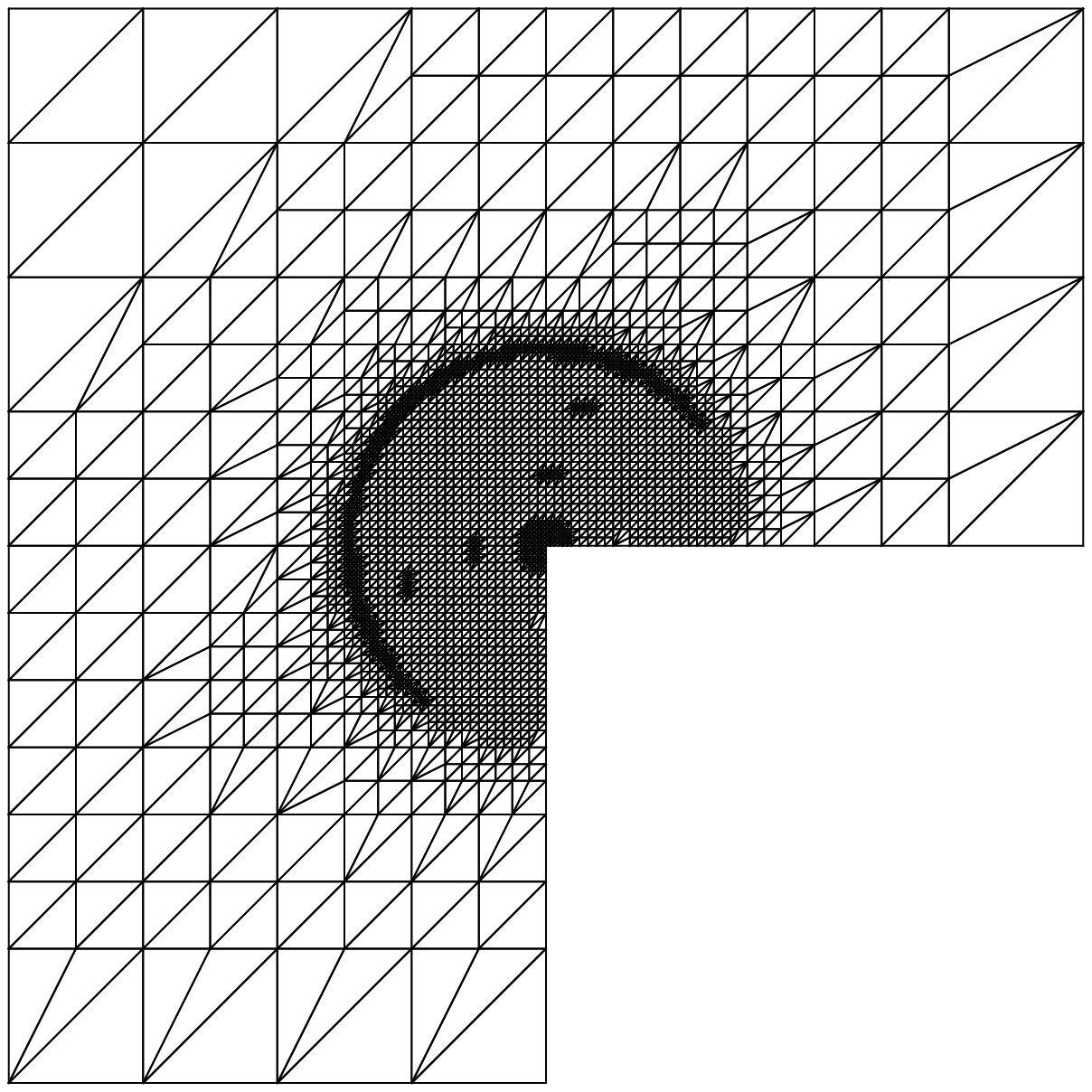}
\end{center}
\caption{Sequence of refined meshes.\label{fig:meshes}}
\end{figure}

\begin{figure}[h]
\begin{center}
\includegraphics[height=9cm]{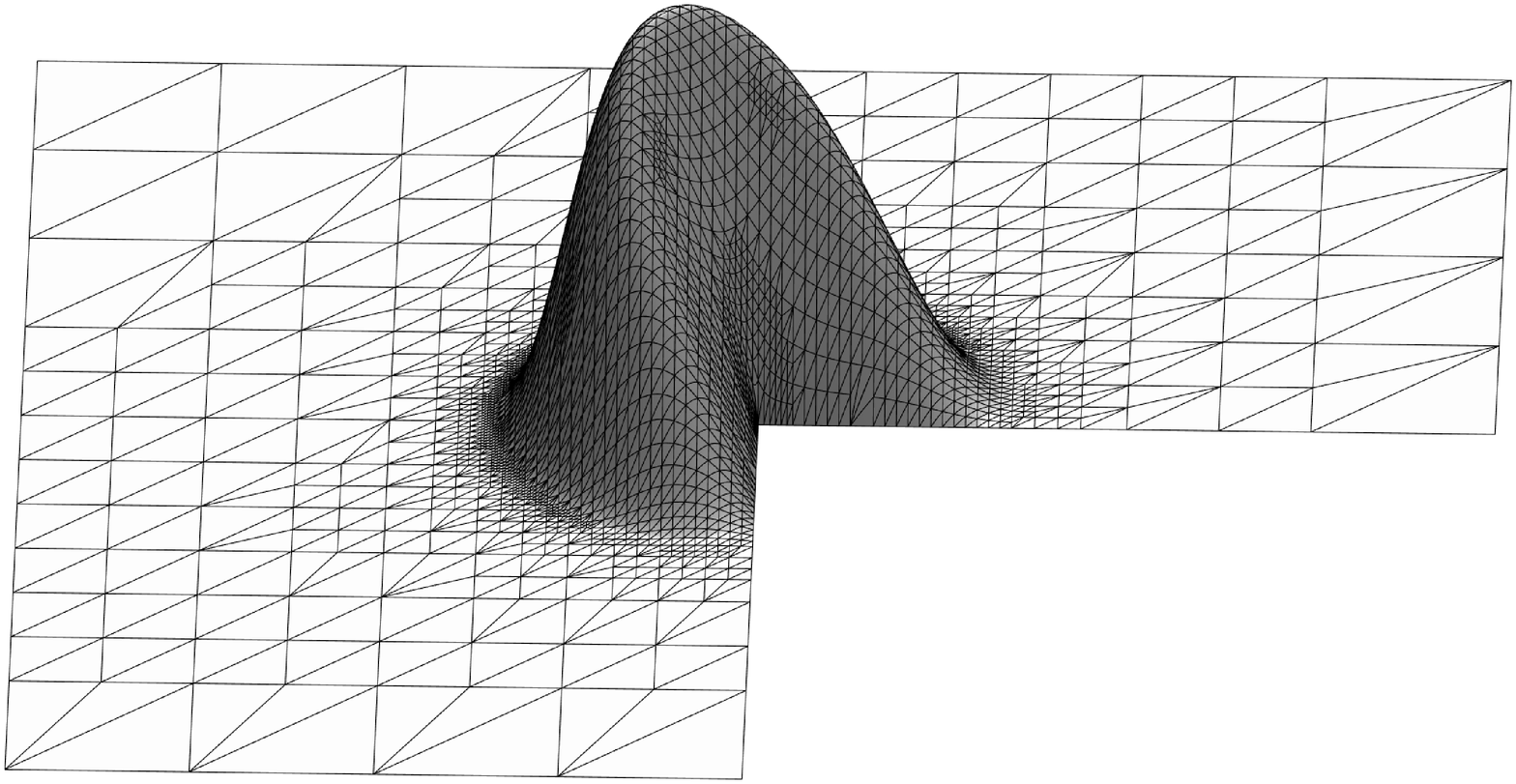}
\end{center}
\caption{Elevation of the discrete solution, nonsmooth case.\label{fig:elnonsmo}}
\end{figure}

For this example we plot, in Fig. \ref{fig:nonsmootherror} the error on consecutive adaptively refined meshes, using the minimum meshsize as a measure of $h$.
We note the suboptimal convergence and that the indicator still approximately follows the $H^1$ error asymptotically. In Fig. \ref{fig:meshes}
we show the corresponding sequence of refined meshes, and in Fig. \ref{fig:elnonsmo} we show an elevation of the approximate solution on the final mesh in the sequance.

\bibliographystyle{plain}
\bibliography{contact}
\newpage

\end{document}